\documentclass[a4paper,11pt,twoside,reqno]{amsart}
\numberwithin{equation}{section}

\usepackage[utf8]{inputenc}
\usepackage[plainpages=false,pdfpagelabels=true]{hyperref}
\usepackage{amssymb,amsthm}
\usepackage[margin=1in]{geometry}
\usepackage[dvipsnames]{xcolor}

\newtheorem{Satz}{Theorem}[section]
\newtheorem{Prop}[Satz]{Proposition}
\newtheorem{Lem}[Satz]{Lemma}

\theoremstyle{definition}

\newtheorem{Bem}[Satz]{Remark}

\newcommand{\tr}{\operatorname{Tr}}
\newcommand{\sff}{\mathrm{I\!I}}

\newcommand{\dv}{~dv_g}
\allowdisplaybreaks[1]

\renewcommand{\epsilon}{\varepsilon}

\newcommand{\R}{\ensuremath{\mathbb{R}}}

\begin{document}

\title{On conservation laws for polyharmonic maps}
\author{Volker Branding}
\address{Oskar-Morgenstern-Platz 1, 1090 Vienna, Austria}
\curraddr{}
\email{volker.branding@univie.ac.at}
\thanks{The author gratefully acknowledge the support of the Austrian Science Fund (FWF) through the project "Geometric Analysis of Biwave Maps" (DOI: 10.55776/P34853).
}
\date{}

\begin{abstract}
This article provides an overview on various conservation laws for polyharmonic maps
between Riemannian manifolds. Besides recalling 
that the variation of the energy for polyharmonic maps with respect to the domain metric gives rise to the stress-energy tensor, we also show how the presence of a Killing vector field on the target manifold leads to a conservation law. For harmonic and biharmonic maps we also
point out a number of applications of such conservation laws.
\end{abstract}

\maketitle

\section{Introduction and Results}
One of the most famous variational problems for maps between two Riemannian manifolds \((M,g)\) and \((N,h)\)
is the one of harmonic maps. For a map \(\phi\colon M\to N\) one defines its energy by
\begin{align}
\label{energy-harmonic}
E(\phi):=\frac{1}{2}\int_M|d\phi|^2\dv.
\end{align}

The critical points of \eqref{energy-harmonic} are characterized by the vanishing 
of the so-called \emph{tension field}, that is
\begin{align}
\label{equation-harmonic}
0=\tau(\phi):=\tr_g\bar\nabla d\phi,
\end{align}
where \(\bar\nabla\) represents the connection on the pull-back bundle \(\phi^\ast TN\).
Solutions of \eqref{equation-harmonic} are called \emph{harmonic maps}.
The harmonic map equation constitutes a second order semilinear elliptic partial differential equation
for which many results have been established, for the current status of research see  \cite{MR2044031} and \cite{MR2389639}.

A fourth order variational problem for maps between Riemannian manifolds that has received growing attention over the last decades is the one of biharmonic maps. To define the latter, one considers the \emph{bienergy} of a map defined as follows
\begin{align}
\label{energy-biharmonic}
E_2(\phi):=\frac{1}{2}\int_M|\tau(\phi)|^2\dv.
\end{align}
The critical points of \eqref{energy-biharmonic} are called \emph{biharmonic maps}
and are determined by the vanishing of the bitension field \(\tau_2(\phi)\),
that is
\begin{align}
\label{bitension}
0=\tau_2(\phi):=-\bar\Delta\tau(\phi)-\tr_g R^N(d\phi(\cdot),\tau(\phi))d\phi(\cdot).
\end{align}
Here, \(R^N\) represents the curvature tensor of the manifold \(N\).
The biharmonic map equation \eqref{bitension} comprises a semilinear elliptic partial differential equation of fourth order and,
due to the higher number of derivatives, additional technical problems arise in its analysis. 

A direct inspection of the biharmonic map equation \eqref{bitension} reveals that every harmonic map is automatically biharmonic. For this reason, one is interested in finding biharmonic maps which are non-harmonic and these are usually referred to as \emph{proper biharmonic}.

In the case that \(M\) is compact and \(N\) has non-positive sectional curvature the maximum principle
shows that all biharmonic maps are necessarily harmonic \cite{MR886529}. For this reason most of the research carried out on biharmonic maps considers the case of a target with positive curvature and many examples of proper biharmonic maps into the sphere have been obtained.

For the current status of research on biharmonic maps between Riemannian manifolds one may consult the recent book \cite{MR4265170} and the survey \cite{MR4410183}.

Besides harmonic and biharmonic maps we will also focus on the following $k$-order versions of the energy functional:

If $k=2s$, $s \geq 1$, we set
\begin{align}
\label{2s-energia}
E_{2s}(\phi)
:=&\frac{1}{2}\int_M |\bar{\Delta}^{s-1}\tau(\phi)|^2\dv.
\end{align}
If $k=2s+1$, $s \geq 1$, we consider
\begin{align}
\label{2s+1-energia}
E_{2s+1}(\phi)
:=&\frac{1}{2}\int_M|\bar\nabla\bar{\Delta}^{s-1}\tau(\phi)|^2 \dv.
\end{align}

A \textit{polyharmonic map of order} $k$ (briefly, a $k$-\textit{harmonic map}) is a critical point of the $k$-energy functional $E_k$. 

The Euler-Lagrange equations of \eqref{2s-energia}, \eqref{2s+1-energia} are  as follows
(with \(\bar\Delta^{-1}=0\)):
\begin{enumerate}
 \item The critical points of \eqref{2s-energia} are those which satisfy
\begin{align}
\label{tension-2s}
0=\tau_{2s}(\phi):=&\bar\Delta^{2s-1}\tau(\phi)-R^N(\bar\Delta^{2s-2}\tau(\phi),d\phi(e_j))d\phi(e_j) \\
\nonumber&-\sum_{\ell=1}^{s-1}\bigg(R^N(\bar\nabla_{e_j}\bar\Delta^{s+\ell-2}\tau(\phi),\bar\Delta^{s-\ell-1}\tau(\phi))d\phi(e_j) \\
\nonumber&\hspace{1cm}-R^N(\bar\Delta^{s+\ell-2}\tau(\phi),\bar\nabla_{e_j}\bar\Delta^{s-\ell-1}\tau(\phi))d\phi(e_j)
\bigg).
\end{align}
\item The critical points of \eqref{2s+1-energia} are determined by
\begin{align}
\label{tension-2s+1}
0=\tau_{2s+1}(\phi):=&\bar\Delta^{2s}\tau(\phi)-R^N(\bar\Delta^{2s-1}\tau(\phi),d\phi(e_j))d\phi(e_j)\\
\nonumber&-\sum_{\ell=1}^{s-1}\bigg(R^N(\bar\nabla_{e_j}\bar\Delta^{s+\ell-1}\tau(\phi),\bar\Delta^{s-\ell-1}\tau(\phi))d\phi(e_j) \\
\nonumber&-R^N(\bar\Delta^{s+\ell-1}\tau(\phi),\bar\nabla_{e_j}\bar\Delta^{s-\ell-1}\tau(\phi))d\phi(e_j)
\bigg) \\
&\nonumber-R^N(\bar\nabla_{e_j}\bar\Delta^{s-1}\tau(\phi),\bar\Delta^{s-1}\tau(\phi))d\phi(e_j).
\end{align}
\end{enumerate}
Here, \(\{e_j\}_{j=1,\ldots,\dim M}\) is a local orthornormal frame field
tangent to \(M\).

For the current status of research on polyharmonic maps we refer to \cite{MR4106647},
polyharmonic hypersurfaces in Riemannian space forms have been investigated in \cite{MR4462636},
this analysis has recently been extended to the pseudo-Riemannian case in \cite{MR4552081}.
Unique continuation properties for polyharmonic maps were obtained in \cite{MR4550874}, a structure theorem for polyharmonic maps from complete non-compact 
Riemannian manifolds was established in \cite{MR4184658}.

Throughout this article we employ the following sign conventions.
The Riemannian curvature tensor is defined by
\[R(X,Y)Z:=[\nabla_X,\nabla_Y]Z-\nabla_{[X,Y]}Z,\]
where \(X,Y,Z\) are vector fields.

For the rough Laplacian on \(\phi^\ast TN\), we shall use the convention
\[\bar\Delta:=-\tr_g(\bar\nabla\bar\nabla-\bar\nabla_\nabla).\]
Hence, the sign of the Laplace operator acting on functions is $\Delta f= -f''$ on $\R$.

Whenever employing local coordinates we will use Latin letters to represent coordinates on \(M\)
and Greek letters for coordinates on the target \(N\).

For vector fields on \(N\) we will use the letters \(X,Y,Z,W\)
and most often \(X\) will refer to a Killing vector field,
while for vector fields on \(M\) we reserve the letters \(U,V\).

Throughout this note we will employ the summation convention, i.e. we
tacitly sum over repeated indices.

This article is organized as follows: In Section 2 we recall a number of results
on conservation laws for harmonic and biharmonic maps which highlight their importance
in differential geometry and geometric analysis.
In Section 3 we first recall the stress-energy tensor for polyharmonic maps
and, in a second step, we show how the existence of a Killing vector field on \(N\)
leads to a conservation law for polyharmonic maps.

\section{Conservation laws for harmonic and biharmonic maps}
In this section we collect several known mathematical results which 
we will need throughout this article. These results will mostly
apply to harmonic and biharmonic maps.

We choose local coordinates \(x^i\) on \(\tilde U\subset M\) and \(y^\alpha\) on \(\tilde V\subset N\)
such that \(\phi(\tilde U)\subset \tilde V\).
Recall that in our sign convention the tension field of a map can be expressed 
in terms of such local coordinates as follows
\begin{align}
\label{tension-field-local}
\tau^\alpha(\phi)=-\Delta\phi^\alpha+\Gamma^\alpha_{\beta\gamma}\langle d\phi^\beta,d\phi^\gamma\rangle,
\end{align}
where 
\begin{align*}
\langle d\phi^\beta,d\phi^\gamma\rangle:=g^{ij}\frac{\partial\phi^\beta}{\partial x^i}\frac{\partial\phi^\gamma}{\partial x^j}
\end{align*}
and \(\Gamma^\alpha_{\beta\gamma}\) represent the Christoffel symbols on the manifold \(N\).

In the following we will frequently make use of 
the Lie derivative of the metric \(h\) on \(N\) with respect to a vector field \(X\). 
This is defined as follows
\begin{align}
(\mathcal{L}_Xh)(Y,Z)=\langle\nabla_YX,Z\rangle+\langle\nabla_ZX,Y\rangle,
\end{align}
where \(X,Y,Z\) are vector fields on \(N\).

A vector field \(X\) is called \emph{Killing vector field} if its Lie derivative
with respect to the metric vanishes, that is
\begin{align}
\label{killing-vector-field}
\mathcal{L}_Xh=0.
\end{align}

As a non-trivial Killing vector field encodes 
information on the symmetries of a Riemannian manifold and as symmetries lead 
to conservation laws via Noether's theorem, it is clear that Killing vector fields
will play an important role within this manuscript.

It directly follows that for a Killing vector field \(X\) we have
\begin{align}
\label{skew-symmetry-killing}
\langle\nabla_YX,Z\rangle+\langle\nabla_ZX,Y\rangle=0.
\end{align}

The next statement relates Killing vector fields to the Riemannian curvature tensor.

\begin{Lem}
Let \(X\) be a Killing vector field on a Riemann manifold \((N,h)\).
Then, the Riemann curvature tensor can be calculated as
\begin{align}
\label{curvature-killing}
\nabla^2_{Y,Z}X=-R^N(X,Y)Z,
\end{align}
where \(Y,Z\) are vector fields on \(N\).
\end{Lem}

\begin{proof}
A proof can be found in \cite[p. 242, Lemma 33]{MR2243772}.
\end{proof}

In the following Lemma we assume that the Lie derivative acts on the metric \(h\) on 
\(N\) and study the behavior of the Christoffel symbols under the action 
of the Lie derivative.

\begin{Lem}
The Lie derivative of the Christoffel symbols with respect to a vector field \(X\)
is given by
\begin{align}
\label{lie-derivative-christoffel}
\mathcal{L}_X\Gamma^\alpha_{\beta\gamma}=\nabla_\beta\nabla_\gamma X^\alpha-R^\alpha_{\beta\gamma\delta}X^\delta,
\end{align}
where \(R^\alpha_{\beta\gamma\delta}\) represent the components
of the Riemann curvature tensor on \(N\).   
\end{Lem}
\begin{proof}
For a proof we refer to \cite[p. 8]{MR88769}.
\end{proof}

Throughout this manuscript we use the notation \(X\circ\phi\)
to highlight that \(X\) is a vector field along the map \(\phi\).

\subsection{The case of harmonic maps}
To obtain a conservation law for harmonic maps we will calculate the Lie derivative of \(|d\phi|^2\) with respect to a vector field \(X\) on \(N\), where the Lie derivative
is supposed to act on the metric \(h\) of \(N\).
More precisely, denoting by \(\{e_i\},i=1,\ldots,m\) a local orthornormal frame tangent to \(M\), we find
\begin{align*}
\mathcal{L}_X|d\phi|^2=&(\mathcal{L}_Xh)
(d\phi(e_i),d\phi(e_i)) \\
=&2\langle d\phi(e_i),\nabla^N_{d\phi(e_i)}X\rangle \\
=&2\langle d\phi(e_i), \bar\nabla_{e_i}(X\circ\phi)\rangle \\
=&2\nabla_{e_i}\langle d\phi(e_i),X\circ\phi\rangle
-2\langle \tau(\phi),X\circ\phi\rangle \\
=&2\operatorname{div} J^1(\phi)-2\langle \tau(\phi),X\circ\phi\rangle,
\end{align*}
where the vector field \(J^1(\phi)\) is defined as follows
\begin{align}
\label{harmonic-current}
J^1(\phi)=\langle d\phi(\cdot), X\circ\phi\rangle^\sharp.
\end{align}
Here, \(\sharp\) represents the musical isomorphism.

Hence, we obtain the following statement which was first obtained in the seminal work of H\'{e}lein \cite{MR1085633}.
\begin{Prop}
\label{proposition-harmonic-conserved-current}
Let \(\phi\colon M\to N\) be a harmonic map. Assume that \(N\) has a Killing vector field \(X\).
Then the vector field \eqref{harmonic-current} is .
\end{Prop}

\begin{proof}
We choose a local geodesic frame field \(\{e_i\},i=1,\ldots,\dim M\)
such that at a fixed point \(p\in M\) we have \(\nabla_{e_j}e_i=0,i,j=1,\ldots,\dim M\).
We calculate
\begin{align*}
\operatorname{div} (J^1(\phi))=\nabla_{e_i}\langle d\phi(e_i),X\circ\phi\rangle=
\langle\underbrace{\tau(\phi)}_{=0},X\circ\phi\rangle
+\langle d\phi(e_i), \bar\nabla_{e_i}(X\circ\phi)\rangle=0,
\end{align*}
and the second term vanishes since \(X\) is a Killing vector field, that is a solution of \eqref{killing-vector-field}.
\end{proof}

In the case of \(u:=\iota\circ\phi\colon M\to\mathbb{S}^n\subset\R^{n+1}\) the vector field \eqref{harmonic-current} 
can be written as
\begin{align}
J^1(u)=u\wedge\nabla u,
\end{align}
where \(\wedge\) represents the exterior product on \(\R^{n+1}\)
and \(\iota\) the canonical embedding of \(\mathbb{S}^n\) into \(\mathbb{R}^{n+1}\).

Note that for \(\mathbb{S}^n\subset\mathbb{R}^{n+1}\) the Killing vector fields
can be represented by antisymmetric matrices, if we restrict these to the sphere
orthogonally, then these will represent Killing vector fields on \(\mathbb{S}^n\).

Moreover, for a spherical target and in our sign convention
the tension field acquires the form
\begin{align}
\label{harmonic-sphere}
d\iota(\tau(\phi))=-\Delta u+|\nabla u|^2u,
\end{align}
where \(\iota\colon\mathbb{S}^n\to\R^{n+1}\).

The next proposition is well-known but underlines the importance of conservation laws
in the analysis of harmonic maps.
\begin{Prop}
Let \(u\colon M\to\mathbb{S}^n\subset\mathbb{R}^{n+1}\) be a smooth map.
Then the following two statements are equivalent:
\begin{enumerate}
\item The map \(u\) is harmonic. 
\item The vector field \(J^1(u)\) is .
\end{enumerate}
\end{Prop}

\begin{proof}
Suppose that \(u\colon M\to\mathbb{S}^n\subset\mathbb{R}^{n+1}\) is a harmonic map, then the vector field 
\(J^1(u)\)
is  due to Proposition \ref{proposition-harmonic-conserved-current}.
For the converse direction suppose that \(J^1(u)\) is  such that
\begin{align*}
u\wedge\Delta u=0.
\end{align*}
This implies that 
\begin{align*}
\Delta u=\lambda u
\end{align*}
for some non-vanishing \(\lambda\in\R\). Taking the scalar product with \(u\)
and using \(\langle u,\Delta u\rangle=|\nabla u|^2\) we find \(\lambda=|\nabla u|^2\).
Consequently, \(u\) solves the equation for harmonic maps in the case of a spherical
target \eqref{harmonic-sphere}.
\end{proof}

\begin{Bem}
The possibility of rewriting the harmonic map equation as a conservation law has led to strong
mathematical results. For example, it is the key ingredient in the existence result
of Shatah for wave maps from Minkowski space to spheres \cite{MR933231}
and the regularity result of H\'{e}lein for harmonic maps from surfaces to spheres \cite{MR1078114}.
\end{Bem}

Let us also make a remark on the relation between conservation laws and integrable systems.
\begin{Bem}
Assume that \(M\) is two-dimensional.
Choosing local coordinates \(x_i,i=1,2\) on an open subset \(\Omega\subset M\)
and again considering the case of a spherical target, we may locally express \(J^1(u)\) as
follows
\begin{align}
\label{current-harmonic-sphere}
(J^1)^{\alpha\beta}_i(u)=u^\alpha_iu^\beta-u^\beta_iu^\alpha,
\end{align}
where we use the notation \(u^\alpha_i:=\frac{\partial u^\alpha}{\partial x^i}\).

By a direct calculation, exploiting the fact that we are considering a spherical target,
we find that \eqref{current-harmonic-sphere} satisfies 
\begin{align}
\label{zero-curvature-harmonic-sphere}
\partial_{x_1}J^1_{x_2}-\partial_{x_2}J^1_{x_1}=2[J^1_{x_1},J^1_{x_2}].
\end{align}
The condition \eqref{zero-curvature-harmonic-sphere} can be identified as a zero curvature equation.
Similar zero curvature conditions can also be found if the target is a symmetric space instead of a sphere.
This construction points out a beautiful connection between harmonic maps to symmetric spaces
and the theory of integrable systems, see for example the book \cite{MR1913803} for more details.
\end{Bem}

For additional results concerning conservation laws for harmonic maps we refer
to \cite{MR716320}, \cite{MR1132860}.

\subsection{The biharmonic case}
In this section we extend the previous analysis to the case of biharmonic maps between Riemannian manifolds.
Some of these results have already been derived in \cite{MR4142862}, \cite{MR2094320}
where a different notation and sign convention was used.

\begin{Lem}
Let \(\phi\colon M\to N\) be a smooth map.
The Lie derivative of \(|\tau(\phi)|^2\) with respect to a vector field \(X\) on \(N\) is given by
\begin{align}
\label{lie-derivative-tension-squared}
\mathcal{L}_X|\tau(\phi)|^2=&(\mathcal{L}_Xh)(\tau(\phi),\tau(\phi))
-2\langle X\circ\phi,\tau_2(\phi)\rangle\\
\nonumber&+2\nabla_{e_j}\big(\langle\bar\nabla_{e_j} (X\circ\phi),\tau(\phi)\rangle
-\langle X\circ\phi,\bar\nabla_{e_j}\tau(\phi)\rangle\big),
\end{align}
where \(\{e_i\},i=1,\ldots,\dim M\) is an orthornormal basis of \(TM\).
\end{Lem}
\begin{proof}
We calculate
\begin{align*}
\mathcal{L}_X|\tau(\phi)|^2=(\mathcal{L}_Xh)(\tau(\phi),\tau(\phi))
+2h(\mathcal{L}_X\tau(\phi),\tau(\phi)).
\end{align*}
Since the tension field depends on the metric \(h\) via the connection of \(N\)
we will get a non-vanishing contribution from \(\mathcal{L}_X\tau(\phi)\).
Using the local expression of the tension field \eqref{tension-field-local} we find
\begin{align*}
\mathcal{L}_X\tau^\alpha(\phi)=&\langle d\phi^\beta,d\phi^\gamma\rangle \mathcal{L}_X\Gamma^\alpha_{\beta\gamma}\\
=&\langle d\phi^\beta,d\phi^\gamma\rangle(\nabla_\beta\nabla_\gamma X^\alpha-R^\alpha_{\beta\gamma\delta}X^\delta),
\end{align*}
where we used \eqref{lie-derivative-christoffel} in the second step.

Hence, we obtain
\begin{align*}
h(\mathcal{L}_X\tau(\phi),\tau(\phi))=&
\langle\bar\nabla_{e_j}\bar\nabla_{e_j}(X\circ\phi),\tau(\phi)\rangle
-\langle R^N(d\phi(e_j),X\circ\phi)d\phi(e_j),\tau(\phi)\rangle \\
=&-\langle\bar\Delta (X\circ\phi),\tau(\phi)\rangle
-\langle R^N(d\phi(e_j),\tau(\phi))d\phi(e_j),X\circ\phi\rangle \\
=&\nabla_{e_j}\big(\langle\bar\nabla_{e_j} (X\circ\phi),\tau(\phi)\rangle-\langle X\circ\phi,\bar\nabla_{e_j}\tau(\phi)\rangle\big) \\
&-\langle X\circ\phi,\bar\Delta\tau(\phi)\rangle
-\langle R^N(d\phi(e_j),\tau(\phi))d\phi(e_j),X\circ\phi\rangle
\end{align*}
completing the proof. 
\end{proof}

Note that if \(X\) is a Killing vector field then the first term on the right hand side in equation \eqref{lie-derivative-tension-squared}
vanishes. Moreover, if we also assume that \(\phi\) is a biharmonic map then \eqref{lie-derivative-tension-squared}
immediately gives us a conservation law as follows:

\begin{Prop}
\label{prop:biharmonic-vf}
Let \(\phi\colon M\to N\) be a smooth biharmonic map. Assume that \(N\) has a Killing vector field \(X\).
Then the following vector field is divergence-free
\begin{align}
\label{biharmonic-vf}
J^2(\phi)=\langle\bar\nabla (X\circ\phi),\tau(\phi)\rangle^\sharp
-\langle    X\circ\phi,\bar\nabla\tau(\phi)\rangle^\sharp,
\end{align}
where \(\sharp\) denotes the musical isomorphism.
\end{Prop}
\begin{proof}
We choose a local geodesic frame field \(\{e_i\},i=1,\ldots,\dim M\)
such that at a fixed point \(p\in M\) we have \(\nabla_{e_j}e_i=0,i,j=1,\ldots,\dim M\).
We calculate
\begin{align*}
\operatorname{div}J^2(\phi)=&\nabla_{e_j}\big(\langle\bar\nabla_{e_j} (X\circ\phi),\tau(\phi)\rangle-\langle X\circ\phi,\bar\nabla_{e_j}\tau(\phi)\rangle\big) \\
=&\langle\bar\nabla_{e_j}\bar\nabla_{e_j} (X\circ\phi),\tau(\phi)\rangle
-\langle X\circ\phi,\bar\nabla_{e_j}\bar\nabla_{e_j}\tau(\phi)\rangle\\
=&-\langle R^N(X\circ\phi,d\phi(e_j))d\phi(e_j),\tau(\phi)\rangle
+\langle X\circ\phi,\bar\Delta\tau(\phi)\rangle \\
=&-\langle X\circ\phi,\tau_2(\phi)\rangle,
\end{align*}
where we have used \eqref{curvature-killing} after the third equals sign,
proving the assertion.
\end{proof}

Now, let us briefly recall the structure of the biharmonic map equation for a spherical target \(\mathbb{S}^n\subset\mathbb{R}^{n+1}\).
Again, recall that in this case the tension field simplifies as
\begin{align*}
d\iota(\tau(\phi))=-\Delta u+|\nabla u|^2u.
\end{align*}

Hence, for \(u:=\iota\circ\phi\colon M\to\mathbb{S}^n\subset\mathbb{R}^{n+1}\), the bienergy acquires the simpler form
\begin{align*}
E_2(u)=\frac{1}{2}\int_M\big(|\Delta u|^2-|\nabla u|^4\big)\dv
\end{align*}
and we now present an extrinsic characterization of its critical points.

\begin{Prop}
Let \(u\colon M\to\mathbb{S}^n\subset\mathbb{R}^{n+1}\) be a smooth map.
Then, the critical points of \(E_2(u)\) are characterized by
\begin{align}
\label{intrinsic-biharmonic-sphere}
\Delta^2u=(-|\Delta u|^2+\Delta|\nabla u|^2+2\langle\nabla u,\nabla\Delta u\rangle-2|\nabla u|^4)u
-2\nabla(|\nabla u|^2\nabla u).
\end{align}
\end{Prop}

\begin{proof}
For a proof one can consult \cite[Lemma 2.1]{MR2094320}.
\end{proof}

In the case of a spherical target the next Lemma shows that the equation for biharmonic maps
\eqref{intrinsic-biharmonic-sphere} is equivalent to a conservation law, this result was first presented in 
\cite[Lemma 2.2]{MR2094320}.
In this setup the vector field \eqref{biharmonic-vf} acquires the form  
\begin{align}
\label{noether-intrinsic-biharmonic-sphere}
J^2(u)=-\nabla\Delta u\wedge u+\Delta u\wedge\nabla u+2|\nabla u|^2\nabla u\wedge u,
\end{align}
where \(\wedge\) represents the exterior product in \(\R^{n+1}\).

\begin{Lem}
Let \(u\colon M\to\mathbb{S}^n\subset\R^{n+1}\) be a smooth map.
Then the following two statements are equivalent
\begin{enumerate}
\item The map \(u\) is a biharmonic map, that is, a solution of \eqref{intrinsic-biharmonic-sphere}. 
\item The vector field \(J^2(u)\) is divergence-free.
\end{enumerate}
\end{Lem}
\begin{proof}
First of all, we note that by taking the wedge product of \eqref{intrinsic-biharmonic-sphere} with \(u\)
we obtain the equation
\begin{align}
\label{divergence-J-biharmonic}
\Delta^2u\wedge u=-2\nabla(|\nabla u|^2\nabla u)\wedge u=-2\nabla(|\nabla u|^2\nabla u\wedge u).
\end{align}
Then, a direct calculation reveals that 
\begin{align*}
\operatorname{div}J^2(u)=\Delta^2 u\wedge u+2\nabla(|\nabla u|^2\nabla u)\wedge u
\end{align*}
showing that the vector field \(J^2(u)\) defined in \eqref{noether-intrinsic-biharmonic-sphere} is divergence-free
whenever we have a solution of \eqref{intrinsic-biharmonic-sphere}.

To show the converse direction we assume that the vector field \(J^2(u)\) is divergence-free
which means that \eqref{divergence-J-biharmonic} holds
and also implies that
\begin{align*}
\Delta^2u+2\nabla(|\nabla u|^2\nabla u)=\lambda u
\end{align*}
for some non-vanishing \(\lambda\in\R\). To determine \(\lambda\) we multiply the above
equation with \(u\) and use the identities
\begin{align*}
\langle\Delta^2u,u\rangle&=\Delta|\nabla u|^2-|\Delta u|^2+2\langle\nabla u,\nabla\Delta u\rangle, \\
\nonumber \langle \nabla(|\nabla u|^2\nabla u),u\rangle &=-|\nabla u|^4
\end{align*}
which follow from differentiating \(|u|^2=1\). The proof is now complete.
\end{proof}

We finish our analysis of \eqref{noether-intrinsic-biharmonic-sphere} by making some comments if one can expect to find a zero-curvature
representation for biharmonic maps to spheres similar to the results for harmonic maps
presented in the previous subsection.

To this end let \(\Omega\subset M\) be a connected domain on which we choose local coordinates \(x^i,1\leq i\leq \dim M\) and consider a spherical target.
In this setup we can express \eqref{noether-intrinsic-biharmonic-sphere} as follows
\begin{align*}
(J^2)_i^{\alpha\beta}(u)=u^\alpha_{jji}u^\beta-u^\beta_{jji}u^\alpha
-u^\alpha_{jj}u^\beta_i+u^\beta_{jj}u^\alpha_i
+2u^\gamma_ju^\gamma_ju^\alpha_iu^\beta
-2u^\gamma_ju^\gamma_ju^\beta_iu^\alpha.
\end{align*}

However, a direct but lengthy calculation shows that, in general,
the above current does not satisfy the equation of a connection with zero-curvature,
i.e.

\begin{align}
\label{local-current-biharmonic}
\partial_{x_i}J^2_{x_j}-\partial_{x_j}J^2_{x_i}\neq
2[J^2_{x_i},J^2_{x_j}],\qquad 1\leq i,j\leq \dim M.
\end{align}

This fact can easily be seen even without performing a detailed calculation.
The order of derivatives on the left hand side of \eqref{local-current-biharmonic} is the order of derivatives
being contained in the current \(J^2\) plus one additional derivative, while the order 
on the right hand side is twice the order of derivatives appearing in the current
as it enters quadratically into the equation.
It is not hard to realize that in the case of harmonic maps the order of derivatives
in the current is one, such that in the equation for a zero-curvature representation
the order on both sides is two.
In the biharmonic case the orders do not match up, on the left hand side it is four
while it is six on the right hand side.

Although the vector field \eqref{local-current-biharmonic} does not help to find a zero curvature representation for biharmonic maps it can be considered as an important ingredient in the regularity theory for fourth order equations developed by Lamm and
Rivi\`{e}re in \cite{MR2398228}.

\subsection{Conservation laws for biharmonic hypersurfaces in spheres}
An important class of biharmonic maps which has received a lot of attention from differential geometers
is the one of biharmonic hypersurfaces in spheres. The current results on biharmonic hypersurfaces available 
all support the conjecture that biharmonic hypersurfaces in spheres must have constant mean curvature (CMC) but
the conjecture in its general form is still open.

In order to extend the general conservation law for biharmonic maps,
given by Proposition \ref{prop:biharmonic-vf},
to biharmonic hypersurfaces in spheres 
let us recall a number of geometric facts on hypersurfaces \(M^m\) in a Riemannian manifold 
\(N^{m+1}\).
The connections on \(M^m\) and \(N^{m+1}\) are related by the following formula
\begin{align}
\label{dfn:sff}
\nabla^N_YZ=\nabla^M_YZ +\sff(Y,Z),
\end{align}
where \(\sff\) represents the second fundamental form of the hypersurface
and \(Y,Z\) are vector fields on \(M\).

Let \(\nu\) be the global unit normal of the hypersurface \(M^m\), then 
its shape operator \(A\) is given by
\begin{align}
\label{dfn:shape}
\nabla^N_{Z}\nu=-A(Z),
\end{align}
where \(Z\) is a vector field on \(M\).

The shape operator \(A\) and the second fundamental form \(\sff\)
are related by
\begin{align}
\label{relation-shape-sff}
\sff(Y,Z)=\langle A(Y),Z\rangle\nu.
\end{align}

If \(\phi\colon M^{m}\to N^{m+1}\) is an isometric immersion
the tension field acquires the form \(\tau(\phi)=mf\nu\),
where \(f\) represents the mean curvature function of the hypersurface.
Recall that we can compute the mean curvature function by \(f=\frac{1}{m}\tr A\).

Also, it is well-known that $\tr_g (\nabla A) (\cdot,\cdot)= m\operatorname{grad} f$,
see for example \cite[Lemma 2.1]{MR4462636}.

A hypersurface in the Euclidean sphere \(\mathbb{S}^{m+1}\) is biharmonic 
if the following system of equations holds
\begin{align}
\label{eqn:hypersurface}
\Delta f=(m-|A|^2)f,\qquad A(\operatorname{grad}f)=-\frac{m}{2}f\operatorname{grad}f,
\end{align}
which arises from splitting the equation for a biharmonic hypersurface into its
normal and its tangential part, see \cite{MR2004799}.

Noting that in the case of a hypersurface
\begin{align*}
\bar\nabla_{U}\tau(\phi)=m U(f)\nu-mf A(U), 
\end{align*}
for all vector fields \(U\) of \(M\),
we find that for a biharmonic hypersurface we have to define
\begin{align}
\label{J-biharmonic-hyper}
J^{2,hyp}(\phi):=mf\langle \bar\nabla (X\circ\phi),\nu\rangle^\sharp
+mf\langle X\circ\phi,A(\cdot)\rangle^\sharp,
\end{align}
where \(\sharp\) denotes the musical isomorphism and \(X\) is a Killing
vector field on \(\mathbb{S}^{m+1}\).

Concerning the vector field \eqref{J-biharmonic-hyper} we prove the following result:
\begin{Satz}
Let \(\phi\colon M^m\to\mathbb{S}^{m+1}\) be a biharmonic hypersurface.
Then, the vector field \(J^{2,hyp}(\phi)\) defined in \eqref{J-biharmonic-hyper} is divergence-free.
\end{Satz}

\begin{proof}
We choose a local geodesic frame field \(\{e_i\},i=1,\ldots,\dim M\)
such that at a fixed point \(p\in M\) we have \(\nabla_{e_j}e_i=0,i,j=1,\ldots,\dim M\).
A direct calculation shows that
\begin{align*}
\operatorname{div}J^{2,hyp}(\phi)=&mf\langle\bar\nabla_{e_i}\bar\nabla_{e_i}(X\circ\phi),\nu\rangle
+2m\langle X\circ\phi,A(\operatorname{grad}f)\rangle 
+mf\langle X\circ\phi,\nabla_{e_i}(A(e_i))\rangle.
\end{align*}

Since we have
\begin{align*}
R^{\mathbb{S}^{m+1}}(Y,Z)W=\langle Z,W\rangle Y-\langle Y,W\rangle Z
\end{align*}
for all vector fields \(Y,Z,W\in\Gamma(T\mathbb{S}^{m+1})\)
we get using \eqref{curvature-killing} that
\begin{align*}
\bar\nabla_{e_i}\bar\nabla_{e_i} (X\circ\phi)=&-R^{\mathbb{S}^{m+1}}((X\circ\phi),d\phi(e_i))d\phi(e_i) \\
=&-|d\phi|^2 X\circ\phi+\langle (X\circ\phi),d\phi(e_i))\rangle d\phi(e_i)
\end{align*}
and we may conclude that
\begin{align*}
\langle\bar\nabla_{e_i}\bar\nabla_{e_i} (X\circ\phi),\nu\rangle=0.
\end{align*}
Moreover, we find
\begin{align*}
\nabla_{e_i}(A(e_i))=&(\nabla A)(e_i,e_i)+\sff(e_i,A(e_i)) \\
=&m\operatorname{grad} f+|A|^2\nu,
\end{align*}
where we used \eqref{relation-shape-sff}. Using these identities we get
\begin{align*}
\operatorname{div}J^{2,hyp}(\phi)=&
2m\langle X\circ\phi, A(\operatorname{grad}f)+\frac{m}{2}f\operatorname{grad}f\rangle
\end{align*}
and thus it is clear that the vector field \(J^{2,hyp}(\phi)\) has vanishing divergence
whenever the second equation of the system \eqref{eqn:hypersurface} holds.
\end{proof}

We would like to point out that the vector field \eqref{J-biharmonic-hyper},
which is divergence-free for a biharmonic hypersurface in \(\mathbb{S}^{m+1}\),
is just a special case of the general statements detailed in
Proposition \ref{prop:biharmonic-vf}.

\section{Conservation laws for polyharmonic maps}
In this section we present two conservation laws that arise in the analysis of polyharmonic maps.
First, we recall the stress-energy tensor for polyharmonic maps, then we will obtain a conservation law
by assuming the the manifold \(N\) has a Killing vector field.
Note that the first conservation law arises from a symmetry on the domain while the second one originates
from a symmetry of the target.

In order to obtain the stress-energy tensor we think of the energies
\eqref{2s-energia}, \eqref{2s+1-energia} for polyharmonic maps as functionals
of the metrics on the domain such that we can vary them with respect to the latter.

\subsection{The stress-energy tensor}
The results in this section have been derived in full detail in \cite{MR4007262},
here we present a short summary of them.

For polyharmonic maps of even order \eqref{tension-2s}, the stress-energy tensor is given by
\begin{align}
\label{even-energy-momentum-tensor}
S_{2s}(U,V):=&g(U,V)\bigg(\frac{1}{2}|\bar\Delta^{s-1}\tau(\phi)|^2-\langle\tau(\phi),\bar\Delta^{2s-2}\tau(\phi)\rangle-\langle d\phi,\bar\nabla\bar\Delta^{2s-2}\tau(\phi)\rangle \\
\nonumber&+\sum_{l=1}^{s-1}(-\langle\bar\Delta^{s-l}\tau(\phi),\bar\Delta^{s+l-2}\tau(\phi)\rangle
+\langle\bar\nabla\bar\Delta^{s-l-1}\tau(\phi),\bar\nabla\bar\Delta^{s+l-2}\tau(\phi)\rangle)
\bigg) \\
\nonumber&-\sum_{l=1}^{s-1}\bigg(\langle\bar\nabla_U\bar\Delta^{s-l-1}\tau(\phi),\bar\nabla_V\bar\Delta^{s+l-2}\tau(\phi)\rangle 
+\langle\bar\nabla_V\bar\Delta^{s-l-1}\tau(\phi),\bar\nabla_U\bar\Delta^{s+l-2}\tau(\phi)\rangle\bigg) \\
\nonumber&+\langle d\phi(U),\bar\nabla_V\bar\Delta^{2s-2}\tau(\phi)\rangle
+\langle d\phi(V),\bar\nabla_U\bar\Delta^{2s-2}\tau(\phi)\rangle,
\end{align}
while for polyharmonic maps of odd order \eqref{tension-2s+1} it is defined as follows
\begin{align}
\label{odd-energy-momentum-tensor}
S_{2s+1}(U,V) 
:=&g(U,V)\bigg(\frac{1}{2}|\bar\nabla\bar\Delta^{s-1}\tau(\phi)|^2
-\langle\tau(\phi),\bar\Delta^{2s-1}\tau(\phi)\rangle
-\langle d\phi,\bar\nabla\bar\Delta^{2s-1}\tau(\phi)\rangle \\
\nonumber &+\sum_{l=1}^{s-1}(-\langle\bar\Delta^{s-l}\tau(\phi),\bar\Delta^{s+l-1}\tau(\phi)\rangle
+\langle\bar\nabla\bar\Delta^{s-l-1}\tau(\phi),\bar\nabla\bar\Delta^{s+l-1}\tau(\phi)\rangle)
\bigg) \\
\nonumber&
-\sum_{l=1}^{s-1}(\langle\bar\nabla_U\bar\Delta^{s-l-1}\tau(\phi),\bar\nabla_V\bar\Delta^{s+l-1}\tau(\phi)\rangle+\langle\bar\nabla_V\bar\Delta^{s-l-1}\tau(\phi),\bar\nabla_U\bar\Delta^{s+l-1}\tau(\phi)\rangle) \\
\nonumber&+\langle d\phi(U),\bar\nabla_V\bar\Delta^{2s-1}\tau(\phi)\rangle
+\langle d\phi(V),\bar\nabla_U\bar\Delta^{2s-1}\tau(\phi)\rangle\\
\nonumber&-\langle\bar\nabla_U\bar\Delta^{s-1}\tau(\phi),\bar\nabla_V\bar\Delta^{s-1}\tau(\phi)\rangle.
\end{align}
Here, \(U\) and \(V\) are vector fields on \(M\). 

It can be directly seen that the stress-energy tensors \(S_{2s}(Y,Z)\) and \(S_{2s+1}(Y,Z)\)
are symmetric. Moreover, they satisfy the following conservation laws:

\begin{Satz}
Let \(\phi\colon M\to N\) be a smooth map.
Then, the stress-energy tensors \eqref{even-energy-momentum-tensor} 
and \eqref{odd-energy-momentum-tensor}
satisfy the following conservation laws
\begin{align*}
\operatorname{div}S_{2s}&=-\langle\tau_{2s}(\phi),d\phi\rangle,\\
\operatorname{div}S_{2s+1}&=-\langle\tau_{2s+1}(\phi),d\phi\rangle.
\end{align*}
In particular, \(S_{2s}\) is divergence-free whenever \(\phi\) is a solution of \eqref{tension-2s}
and \(S_{2s+1}\) is divergence-free whenever \(\phi\) is a solution of \eqref{tension-2s+1}.
\end{Satz}
\begin{proof}
For the precise details see the proofs of Propositions 2.6 and 2.12 in \cite{MR4007262}.
\end{proof}

In Subsection 2.3 of \cite{MR4007262} it is shown how the fact that the stress-energy tensors for polyharmonic maps
are conserved arises from the invariance of the energy functionals \eqref{2s-energia}, \eqref{2s+1-energia}
under diffeomorphisms on the domain and is thus a consequence of Noether's theorem.

A possible application of the stress-energy tensors considered above was pointed out in Section 4 of \cite{MR4007262},
where a classification result for finite energy solutions of the triharmonic map equation (\eqref{tension-2s+1} with \(s=1\)) from Euclidean space was obtained.

\subsection{A conservation law arising from Killing vector fields}
In this subsection we show how the presence of a Killing vector field on the target manifold \(N\)
gives rise to a conservation law for polyharmonic maps. Again, we need to distinguish between the case
of a polyharmonic map of even order and a polyharmonic map of odd order.

Recall that in the case of biharmonic maps we were able to obtain a conservation law by calculating
the Lie derivative of \(|\tau(\phi)|^2\) with respect to a vector field \(X\) on \(N\),
where the Lie derivative is supposed to act on the metric \(h\) on \(N\).
However, it turns out that this approach is no longer feasible for \(k\)-harmonic maps of order \(k\geq 3\)
from a computational point of view.

Luckily, there is also a second constructive approach to derive a conservation law for polyharmonic 
maps inspired by the seminal work of H\'{e}lein for harmonic maps \cite{MR1078114}.

To this end, let \(\eta\) be a smooth function on \(M\) with compact support, that is \(\eta\in C_0^\infty(M)\)
and \(X\) a Killing vector field on \(N\). First, we focus on the case of polyharmonic   maps of even order.

We test \eqref{tension-2s} with \(\eta (X\circ\phi)\) and use integration by parts to manipulate
all contributing terms
\begin{align*}
\int_M\langle\bar\Delta^{2s-1}\tau(\phi),X\circ\phi\rangle\eta\dv=&
\int_M\langle\bar\nabla\bar\Delta^{2s-2}\tau(\phi),\bar\nabla (X\circ\phi)\rangle\eta\dv \\
&+\int_M\langle\bar\nabla\bar\Delta^{2s-2}\tau(\phi),X\circ\phi\rangle\nabla\eta\dv.
\end{align*}
Moreover, we find
\begin{align*}
-\int_M\langle& R^N(\bar\Delta^{2s-2}\tau(\phi),d\phi(e_j))d\phi(e_j),X\circ\phi\rangle\eta\dv \\
=&\int_M\langle\bar\nabla_{e_j}\bar\nabla_{e_j}(X\circ\phi),\bar\Delta^{2s-2}\tau(\phi)\rangle\eta\dv \\
=&-\int_M\langle\bar\nabla (X\circ\phi),\bar\nabla\bar\Delta^{2s-2}\tau(\phi)\rangle\eta\dv 
-\int_M\langle\bar\nabla (X\circ\phi),\bar\Delta^{2s-2}\tau(\phi)\rangle\nabla\eta\dv,
\end{align*}
where we used \eqref{curvature-killing} in the first step.

Now, we manipulate
\begin{align*}
\int_M\langle &R^N(\bar\nabla_{e_j}\bar\Delta^{s+l-2}\tau(\phi),\bar\Delta^{s-l-1}\tau(\phi))d\phi(e_j),X\circ\phi\rangle\eta\dv \\
=&\int_M\langle\bar\nabla_{e_j}\big(\bar\nabla_{\bar\nabla_{e_j}\bar\Delta^{s+l-2}\tau(\phi)}(X\circ\phi)\big),\bar\Delta^{s-l-1}\tau(\phi)\rangle\eta\dv\\
=&-\int_M\langle\bar\nabla_{\bar\nabla\bar\Delta^{s+l-2}\tau(\phi)}(X\circ\phi),\bar\nabla\bar\Delta^{s-l-1}\tau(\phi)\rangle\eta\dv \\
&-\int_M\langle\bar\nabla_{\bar\nabla\bar\Delta^{s+l-2}\tau(\phi)}(X\circ\phi),\bar\Delta^{s-l-1}\tau(\phi)\rangle\nabla\eta\dv.
\end{align*}
Similarly, we find
\begin{align*}
-\int_M\langle &R^N(\bar\Delta^{s+l-2}\tau(\phi),\bar\nabla_{e_j}\bar\Delta^{s-l-1}\tau(\phi))d\phi(e_j),X\circ\phi\rangle\eta\dv \\
=&\int_M\langle\bar\nabla_{e_j}\big(\bar\nabla_{\bar\nabla_{e_j}\bar\Delta^{s-l-1}\tau(\phi)}(X\circ\phi)\big),\bar\Delta^{s+l-2}\tau(\phi)\rangle\eta\dv\\
=&-\int_M\langle\bar\nabla_{\bar\nabla\bar\Delta^{s-l-1}\tau(\phi)}(X\circ\phi),\bar\nabla\bar\Delta^{s+l-2}\tau(\phi)\rangle\eta\dv \\
&-\int_M\langle\bar\nabla_{\bar\nabla\bar\Delta^{s-l-1}\tau(\phi)}(X\circ\phi),\bar\Delta^{s+l-2}\tau(\phi)\rangle\nabla\eta\dv.
\end{align*}
Adding up the different contributions and using the skew-symmetry of the Killing vector field \eqref{skew-symmetry-killing} we find
\begin{align*}
\int_M&\bigg(\langle\bar\nabla_{e_j}\bar\Delta^{2s-2}\tau(\phi),X\circ\phi\rangle
-\langle\bar\nabla_{e_j}(X\circ\phi),\bar\Delta^{2s-2}\tau(\phi)\rangle \\
&+\sum_{l=1}^{s-1}\big(\langle\bar\nabla_{\bar\nabla_{e_j}\bar\Delta^{s+l-2}\tau(\phi)}(X\circ\phi),\bar\Delta^{s-l-1}\tau(\phi)\rangle \\
&+\langle\bar\nabla_{\bar\nabla_{e_j}\bar\Delta^{s-l-1}\tau(\phi)}(X\circ\phi),\bar\Delta^{s+l-2}\tau(\phi)\rangle)\bigg)\nabla_{e_j}\eta\dv \\
&=\int_M\langle J^{2s}(\phi),\nabla\eta\rangle\dv.
\end{align*}

Hence, we can deduce that the vector field, which should be divergence-free for a polyharmonic map of even order, is precisely the following one

\begin{align}
\label{current-polyharmonic-even}
J^{2s}(\phi):=&\langle\bar\nabla\bar\Delta^{2s-2}\tau(\phi),X\circ\phi\rangle
-\langle\bar\nabla (X\circ\phi),\bar\Delta^{2s-2}\tau(\phi)\rangle \\
\nonumber
&+\sum_{l=1}^{s-1}\big(\langle\bar\nabla_{\bar\nabla\bar\Delta^{s+l-2}\tau(\phi)}(X\circ\phi),\bar\Delta^{s-l-1}\tau(\phi)\rangle \\
&\nonumber
+\langle\bar\nabla_{\bar\nabla\bar\Delta^{s-l-1}\tau(\phi)}(X\circ\phi),\bar\Delta^{s+l-2}\tau(\phi)\rangle\big).
\end{align}

This is demonstrated by the following
\begin{Satz}
Let \(\phi\colon M\to N\) be a smooth map and assume that the manifold \(N\) has a Killing vector field \(X\).
Then, the vector field \(J^{2s}(\phi)\) defined in \eqref{current-polyharmonic-even} satisfies 
\begin{align*}
\operatorname{div} J^{2s}(\phi)=-\langle \tau_{2s}(\phi),X\circ\phi\rangle.
\end{align*}
In particular, \(J^{2s}(\phi)\) is divergence-free whenever \(\phi\colon M\to N\) is a polyharmonic map
of even order.
\end{Satz}

\begin{proof}
We choose a local geodesic frame field \(\{e_i\},i=1,\ldots,\dim M\)
such that at a fixed point \(p\in M\) we have \(\nabla_{e_j}e_i=0,i,j=1,\ldots,\dim M\).
A direct calculation shows that
\begin{align*}
\bar\nabla_{e_i}\big(&\langle\bar\nabla_{e_i}\bar\Delta^{2s-2}\tau(\phi),X\circ\phi\rangle
-\langle\bar\nabla_{e_i}(X\circ\phi),\bar\Delta^{2s-2}\tau(\phi)\rangle\big) \\
&=-\langle\bar\Delta^{2s-1}\tau(\phi),X\circ\phi\rangle
+\langle\bar\Delta (X\circ\phi),\bar\Delta^{2s-2}\tau(\phi)\rangle \\
&=-\langle\bar\Delta^{2s-1}\tau(\phi),X\circ\phi\rangle
+\langle R^N(X\circ\phi,d\phi(e_i))d\phi(e_i),\bar\Delta^{2s-2}\tau(\phi)\rangle,
\end{align*}
where we used \eqref{curvature-killing} in the second step.
Moreover, we find
\begin{align*}
\bar\nabla_{e_i}&\big(
\langle\bar\nabla_{\bar\nabla_{e_i}\bar\Delta^{s+l-2}\tau(\phi)}(X\circ\phi),\bar\Delta^{s-l-1}\tau(\phi)\rangle 
+\langle\bar\nabla_{\bar\nabla_{e_i}\bar\Delta^{s-l-1}\tau(\phi)}(X\circ\phi),\bar\Delta^{s+l-2}\tau(\phi)\rangle
\big) \\
=&\langle\bar\nabla_{e_i}\big(\bar\nabla_{\bar\nabla_{e_i}\bar\Delta^{s+l-2}\tau(\phi)}(X\circ\phi)\big),\bar\Delta^{s-l-1}\tau(\phi)\rangle \\
&+\langle\bar\nabla_{\bar\nabla_{e_i}\bar\Delta^{s+l-2}\tau(\phi)}(X\circ\phi),\bar\nabla_{e_i}\bar\Delta^{s-l-1}\tau(\phi)\rangle \\
&+\langle\bar\nabla_{e_i}\big(\bar\nabla_{\bar\nabla_{e_i}\bar\Delta^{s-l-1}\tau(\phi)}(X\circ\phi)\big),\bar\Delta^{s+l-2}\tau(\phi)\rangle\\
&+\langle\bar\nabla_{\bar\nabla_{e_i}\bar\Delta^{s-l-1}\tau(\phi)}(X\circ\phi),\bar\nabla_{e_i}\bar\Delta^{s+l-2}\tau(\phi)\rangle \\
=&-\langle R^N(X\circ\phi,d\phi(e_i))\bar\nabla_{e_i}\bar\Delta^{s+l-2}\tau(\phi),\bar\Delta^{s-l-1}\tau(\phi)\rangle \\
&-\langle R^N(X\circ\phi,d\phi(e_i))\bar\nabla_{e_i}\bar\Delta^{s-l-1}\tau(\phi),\bar\Delta^{s+l-2}\tau(\phi)\rangle,
\end{align*}
where we used both \eqref{skew-symmetry-killing} and \eqref{curvature-killing}.
The claim now follows by using the definition of \(J^{2s}(\phi)\) and the Euler-Lagrange equation
for polyharmonic maps of even order.
\end{proof}

In the case of polyharmonic maps of odd order we also have a conserved vector field, given by
\begin{align}
\label{dfn:j-odd}
J^{2s+1}(\phi):=&\langle\bar\nabla\bar\Delta^{2s-1}\tau(\phi),X\circ\phi\rangle
-\langle\bar\nabla (X\circ\phi),\bar\Delta^{2s-1}\tau(\phi)\rangle \\
\nonumber&+\sum_{l=1}^{s-1}\big(-\langle\bar\nabla_{\bar\nabla\bar\Delta^{s-1}\tau(\phi)}(X\circ\phi),\bar\Delta^{s-1}\tau(\phi)\rangle \\
\nonumber&+\langle\bar\nabla_{\bar\nabla\bar\Delta^{s+l-1}\tau(\phi)}(X\circ\phi),\bar\Delta^{s-l-1}\tau(\phi)\rangle \\
\nonumber&+\langle\bar\nabla_{\bar\nabla\bar\Delta^{s-l-1}\tau(\phi)}(X\circ\phi),\bar\Delta^{s+l-1}\tau(\phi)\rangle
\big).
\end{align}

We can construct the vector field \eqref{dfn:j-odd} using the same arguments as in the even case
such that we do not give the precise derivation here.
Nonetheless, the vector field \eqref{dfn:j-odd} satisfies the following conservation law.

\begin{Satz}
Let \(\phi\colon M\to N\) be a smooth map and assume that the manifold N has a Killing vector field \(X\).
Then, the vector field \(J^{2s+1}(\phi)\) defined in \eqref{dfn:j-odd} satisfies 
\begin{align*}
\operatorname{div} J^{2s+1}(\phi)=-\langle \tau_{2s+1}(\phi),X\circ\phi\rangle.
\end{align*}
In particular, \(J^{2s+1}(\phi)\) is divergence-free whenever \(\phi\colon M\to N\) is a polyharmonic map
of odd order, that is a solution of \eqref{tension-2s}.
\end{Satz}

\begin{proof}
We choose a local geodesic frame field \(\{e_i\},i=1,\ldots,\dim M\)
such that at a fixed point \(p\in M\) we have \(\nabla_{e_j}e_i=0,i,j=1,\ldots,\dim M\).
We calculate
\begin{align*}
\bar\nabla_{e_i}\big(&\langle\bar\nabla_{e_i}\bar\Delta^{2s-1}\tau(\phi),X\circ\phi\rangle
\nonumber-\langle\bar\nabla_{e_i}(X\circ\phi),\bar\Delta^{2s-1}\tau(\phi)\rangle\big)\\
&=-\langle\bar\Delta^{2s}\tau(\phi),X\circ\phi\rangle
+\langle\bar\Delta (X\circ\phi),\bar\Delta^{2s-1}\tau(\phi)\rangle \\
&=-\langle\bar\Delta^{2s}\tau(\phi),X\circ\phi\rangle
+\langle R^N(X\circ\phi,d\phi(e_i))d\phi(e_i),\bar\Delta^{2s-1}\tau(\phi)\rangle.
\end{align*}
In addition, we find
\begin{align*}
\bar\nabla_{e_i}\big(&-\langle\bar\nabla_{\bar\nabla_{e_i}\bar\Delta^{s-1}\tau(\phi)}(X\circ\phi),\bar\Delta^{s-1}\tau(\phi)\rangle
+\langle\bar\nabla_{\bar\nabla_{e_i}\bar\Delta^{s+l-1}\tau(\phi)}(X\circ\phi),\bar\Delta^{s-l-1}\tau(\phi)\rangle \\
&+\langle\bar\nabla_{\bar\nabla_{e_i}\bar\Delta^{s-l-1}\tau(\phi)}(X\circ\phi),\bar\Delta^{s+l-1}\tau(\phi)\rangle
\big)\\
=&-\langle\bar\nabla_{e_i}\big(\bar\nabla_{\bar\nabla_{e_i}\bar\Delta^{s-1}\tau(\phi)}(X\circ\phi)\big),\bar\Delta^{s-1}\tau(\phi)\rangle \\
&-\underbrace{\langle\bar\nabla_{\bar\nabla_{e_i}\bar\Delta^{s-1}\tau(\phi)}(X\circ\phi),\bar\nabla_{e_i}\bar\Delta^{s-1}\tau(\phi)\rangle}_{=0} \\
&+\langle\bar\nabla_{e_i}\big(\bar\nabla_{\bar\nabla_{e_i}\bar\Delta^{s+l-1}\tau(\phi)}(X\circ\phi)\big),\bar\Delta^{s-l-1}\tau(\phi)\rangle \\
&+\langle\bar\nabla_{\bar\nabla_{e_i}\bar\Delta^{s+l-1}\tau(\phi)}(X\circ\phi),\bar\nabla_{e_i}\bar\Delta^{s-l-1}\tau(\phi)\rangle \\
&+\langle\bar\nabla_{e_i}\big(\bar\nabla_{\bar\nabla_{e_i}\bar\Delta^{s-l-1}\tau(\phi)}(X\circ\phi)\big),\bar\Delta^{s+l-1}\tau(\phi)\rangle\\
&+\langle\bar\nabla_{\bar\nabla_{e_i}\bar\Delta^{s-l-1}\tau(\phi)}(X\circ\phi),\bar\nabla_{e_i}\bar\Delta^{s+l-1}\tau(\phi)\rangle \\
=&\langle R^N(X\circ\phi,d\phi(e_i))\bar\nabla_{e_i}\bar\Delta^{s-1}\tau(\phi),\bar\Delta^{s-1}\tau(\phi)\rangle \\
&-\langle R^N(X\circ\phi,d\phi(e_i))\bar\nabla_{e_i}\bar\Delta^{s+l-1}\tau(\phi),\bar\Delta^{s-l-1}\tau(\phi)\rangle 
\\
&-\langle R^N(X\circ\phi,d\phi(e_i))\bar\nabla_{e_i}\bar\Delta^{s-l-1}\tau(\phi),\bar\Delta^{s+l-1}\tau(\phi)\rangle.
\end{align*}
The claim now follows from using the symmetries of the Riemannian curvature tensor and the concrete form
of the equation for polyharmonic maps of odd order \eqref{tension-2s+1}.
\end{proof}

\begin{Bem}
We want to point out that conserved quantities such as \eqref{current-polyharmonic-even} and \eqref{dfn:j-odd}
turn out to be very useful in order to develop a general regularity theory for polyharmonic maps, see \cite{MR4277329} and \cite{MR4279397}. In order to establish such analytic results one only needs a local quantity that is conserved,
hence one does not need to demand the existence of a Killing vector field on \(N\) but can use the Coulomb gauge construction to locally express the equation for polyharmonic maps as a conservation law.
 \end{Bem}

\bibliographystyle{plain}
\bibliography{mybib}

\begin{thebibliography}{10}

\bibitem{MR716320}
Paul Baird.
\newblock {\em Harmonic maps with symmetry, harmonic morphisms and deformations
  of metrics}, volume~87 of {\em Research Notes in Mathematics}.
\newblock Pitman (Advanced Publishing Program), Boston, MA, 1983.

\bibitem{MR1132860}
Paul Baird and Andrea Ratto.
\newblock Conservation laws, equivariant harmonic maps and harmonic morphisms.
\newblock {\em Proc. London Math. Soc. (3)}, 64(1):197--224, 1992.

\bibitem{MR2044031}
Paul Baird and John~C. Wood.
\newblock {\em Harmonic morphisms between {R}iemannian manifolds}, volume~29 of
  {\em London Mathematical Society Monographs. New Series}.
\newblock The Clarendon Press, Oxford University Press, Oxford, 2003.

\bibitem{MR4106647}
V.~Branding, S.~Montaldo, C.~Oniciuc, and A.~Ratto.
\newblock Higher order energy functionals.
\newblock {\em Adv. Math.}, 370:107236, 60, 2020.

\bibitem{MR4552081}
V.~Branding, S.~Montaldo, C.~Oniciuc, and A.~Ratto.
\newblock Polyharmonic hypersurfaces into pseudo-{R}iemannian space forms.
\newblock {\em Ann. Mat. Pura Appl. (4)}, 202(2):877--899, 2023.

\bibitem{MR4142862}
Volker Branding.
\newblock Some analytic results on interpolating sesqui-harmonic maps.
\newblock {\em Ann. Mat. Pura Appl. (4)}, 199(5):2039--2059, 2020.

\bibitem{MR4007262}
Volker Branding.
\newblock The stress-energy tensor for polyharmonic maps.
\newblock {\em Nonlinear Anal.}, 190:111616, 17, 2020.

\bibitem{MR4184658}
Volker Branding.
\newblock A structure theorem for polyharmonic maps between {R}iemannian
  manifolds.
\newblock {\em J. Differential Equations}, 273:14--39, 2021.

\bibitem{MR4550874}
Volker Branding, Stefano Montaldo, Cezar Oniciuc, and Andrea Ratto.
\newblock Unique continuation properties for polyharmonic maps between
  {R}iemannian manifolds.
\newblock {\em Canad. J. Math.}, 75(1):1--28, 2023.

\bibitem{MR4279397}
Fr\'{e}d\'{e}ric~Louis de~Longueville and Andreas Gastel.
\newblock Conservation laws for even order systems of polyharmonic map type.
\newblock {\em Calc. Var. Partial Differential Equations}, 60(4):Paper No. 138,
  18, 2021.

\bibitem{MR4410183}
Dorel Fetcu and Cezar Oniciuc.
\newblock Biharmonic and biconservative hypersurfaces in space forms.
\newblock In {\em Differential geometry and global analysis---in honor of
  {T}adashi {N}agano}, volume 777 of {\em Contemp. Math.}, pages 65--90. Amer.
  Math. Soc., [Providence], RI, [2022] \copyright 2022.

\bibitem{MR1078114}
Fr\'{e}d\'{e}ric H\'{e}lein.
\newblock R\'{e}gularit\'{e} des applications faiblement harmoniques entre une
  surface et une sph\`ere.
\newblock {\em C. R. Acad. Sci. Paris S\'{e}r. I Math.}, 311(9):519--524, 1990.

\bibitem{MR1085633}
Fr\'{e}d\'{e}ric H\'{e}lein.
\newblock Regularity of weakly harmonic maps from a surface into a manifold
  with symmetries.
\newblock {\em Manuscripta Math.}, 70(2):203--218, 1991.

\bibitem{MR1913803}
Fr\'{e}d\'{e}ric H\'{e}lein.
\newblock {\em Harmonic maps, conservation laws and moving frames}, volume 150
  of {\em Cambridge Tracts in Mathematics}.
\newblock Cambridge University Press, Cambridge, second edition, 2002.
\newblock Translated from the 1996 French original, With a foreword by James
  Eells.

\bibitem{MR2389639}
Fr\'{e}d\'{e}ric H\'{e}lein and John~C. Wood.
\newblock Harmonic maps.
\newblock In {\em Handbook of global analysis}, pages 417--491, 1213. Elsevier
  Sci. B. V., Amsterdam, 2008.

\bibitem{MR4277329}
Jasmin H\"{o}rter and Tobias Lamm.
\newblock Conservation laws for even order elliptic systems in the critical
  dimension---a new approach.
\newblock {\em Calc. Var. Partial Differential Equations}, 60(4):Paper No. 125,
  23, 2021.

\bibitem{MR886529}
Guo~Ying Jiang.
\newblock {$2$}-harmonic maps and their first and second variational formulas.
\newblock {\em Chinese Ann. Math. Ser. A}, 7(4):389--402, 1986.
\newblock An English summary appears in Chinese Ann. Math. Ser. B {{\bf{7}}}
  (1986), no. 4, 523.

\bibitem{MR2398228}
Tobias Lamm and Tristan Rivi\`ere.
\newblock Conservation laws for fourth order systems in four dimensions.
\newblock {\em Comm. Partial Differential Equations}, 33(1-3):245--262, 2008.

\bibitem{MR4462636}
Stefano Montaldo, Cezar Oniciuc, and Andrea Ratto.
\newblock Polyharmonic hypersurfaces into space forms.
\newblock {\em Israel J. Math.}, 249(1):343--374, 2022.

\bibitem{MR2004799}
C.~Oniciuc.
\newblock Biharmonic maps between {R}iemannian manifolds.
\newblock {\em An. \c{S}tiin\c{t}. Univ. Al. I. Cuza Ia\c{s}i. Mat. (N.S.)},
  48(2):237--248 (2003), 2002.

\bibitem{MR4265170}
Ye-Lin Ou and Bang-Yen Chen.
\newblock {\em Biharmonic submanifolds and biharmonic maps in {R}iemannian
  geometry}.
\newblock World Scientific Publishing Co. Pte. Ltd., Hackensack, NJ, [2020]
  \copyright 2020.

\bibitem{MR2243772}
Peter Petersen.
\newblock {\em Riemannian geometry}, volume 171 of {\em Graduate Texts in
  Mathematics}.
\newblock Springer, New York, second edition, 2006.

\bibitem{MR933231}
Jalal Shatah.
\newblock Weak solutions and development of singularities of the {${\rm
  SU}(2)$} {$\sigma$}-model.
\newblock {\em Comm. Pure Appl. Math.}, 41(4):459--469, 1988.

\bibitem{MR2094320}
Changyou Wang.
\newblock Remarks on biharmonic maps into spheres.
\newblock {\em Calc. Var. Partial Differential Equations}, 21(3):221--242,
  2004.

\bibitem{MR88769}
Kentaro Yano.
\newblock {\em The theory of {L}ie derivatives and its applications}.
\newblock North-Holland Publishing Co., Amsterdam; P. Noordhoff Ltd., Groningen
  Interscience Publishers Inc., New York, 1957.

\end{thebibliography}

\end{document}